\documentclass
[12pt,a4paper]{article}
\usepackage[cp1251]{inputenc}
\usepackage{amssymb, amsfonts, amsmath, latexsym}
\usepackage[english]{babel}

\begin{document}


\centerline{\bf  Varieties of coarse spaces }\vspace{6 mm}

\normalsize\centerline{\bf  Igor Protasov }\vspace{6 mm}

{\bf Abstract.} A class $\mathfrak{M}$ of  coarse spaces is called a variety if  $\mathfrak{M}$ is closed under formation of subspaces, coarse images and products.  We classify the varieties of coarse spaces and, in particular, show that if a variety  $\mathfrak{M}$ contains an unbounded metric space then  $\mathfrak{M}$ is the variety of all coarse spaces.

\vspace{6 mm}

MSC: 54E35, 08B85.

\vspace{3 mm}

Keywords: coarse structure,coarse space, ballean, varieties of coarse spaces.
\vspace{6 mm}

\section{Introduction}

Following \cite{b10}, we say that a family $\mathcal{E}$ of subsets of $X\times X$ is a {\it coarse structure} on a set $X$ if

\begin{itemize}
\item{} each $\varepsilon \in \mathcal{E}$ contains the diagonal $\vartriangle _{X}$, $\vartriangle _{X}= \{(x,x) : x \in X\}$ ; \vskip 5pt

\item{}  if $\varepsilon, \delta\in\mathcal{E}$ then $\varepsilon \circ\delta\in\mathcal{E}$  and $\varepsilon^{-1}\in\mathcal{E}$ where $\varepsilon \circ\delta = \{(x, y): \exists z ((x,z)\in\varepsilon, (z,y)\in\delta)\}, $  $ \ \varepsilon^{-1}= \{(y,x): (x,y)\in\varepsilon\}$;

\item{}  if $\varepsilon\in\mathcal{E}$ and $\bigtriangleup_{X}\subseteq \varepsilon^{\prime}\subseteq\varepsilon$ then $\varepsilon^{\prime}\in\mathcal{E}$.

\end{itemize}

Each $\varepsilon\in\mathcal{E}$ is called an {\it entourage} of the diagonal.
A subset $\mathcal{E}^{\prime}\subseteq\mathcal{E}$ is called a {\it base} for $\mathcal{E}$ if, for
every $\varepsilon\in\mathcal{E}$ there exists $\varepsilon^{\prime}\in\mathcal{E}^{\prime}$ such that $\varepsilon\subseteq\varepsilon^{\prime}$.

The pair $(X, \mathcal{E})$ is called a {\it coarse space}. For
$x\in X$  and $\varepsilon\in\mathcal{E}$, we denote
 $B(x, \varepsilon)= \{y\in X: (x,y)\in\varepsilon \}$
and say that
 $B(x,\varepsilon)$ is a {\it ball of radius  $\varepsilon$ around $x$.}
 We note that a coarse space can be considered as
 an asymptotic counterpart of a uniform topological space and could
  be defined in terms of balls, see \cite{b7}, \cite{b9}.
 In this case a coarse space is called a {\it ballean}.

A coarse space $(X,\mathcal{E})$ is called {\it connected} if, for any $x, y \in X$,
 there exists  $\varepsilon\in\mathcal{E}$ such that $y\in B(x,\varepsilon)$.
 A subset $Y$ of $X$ is called {\it bounded} if there exist $x\in X$ and $\varepsilon\in\mathcal{E}$
   such that $Y\subseteq B(x, \varepsilon)$.  The coarse
   structure
    $\mathcal{E}=\{\varepsilon\in X\times X: \bigtriangleup_{X}\subseteq\varepsilon\}$ is
     the unique coarse structure such that  $(X,\mathcal{E})$ is connected and bounded.

In what follows, all coarse spaces under consideration are supposed to be {\bf connected}.

Given  a coarse space $(X, \mathcal{E})$, each subset $Y \subseteq X$ has the natural coarse structure
$\mathcal{E}|_{Y}= \{\varepsilon\cap(Y\times Y): \varepsilon\in\mathcal{E} \}$, $(Y, \mathcal{E}|_{Y})$
 is called a {\it subspace} of $(X, \mathcal{E})$.
A subset $Y$ of $X$  is called {\it large} (or {\it  coarsely dense})  if there exists $\varepsilon\in \mathcal{E}$  such that $X= B(Y, \varepsilon)$ where $B(Y, \varepsilon)=\cup_{y\in Y} B(Y, \varepsilon)$.

Let  $(X, \mathcal{E})$, $(X^{\prime}, \mathcal{E}^{\prime})$ be coarse spaces. A mapping
$f: X\longrightarrow X^{\prime}$
is called {\it coarse}  (or {\it  bornologous}  in terminology of \cite{b10}) if, for every
$\varepsilon\in\mathcal{E}$
there exists $\varepsilon^{\prime}\in\mathcal{E}^{\prime}$ such that, for every $x\in X$,  we have
$f(B(x,\varepsilon))\subseteq (B(f(x),\varepsilon^{\prime}))$.
If $f$ is surjective and coarse then $(X^{\prime}, \mathcal{E}^{\prime})$
 is called a {\it coarse image} of $(X, \mathcal{E})$.
If $f$ is a  bijection such that $f$  and $f^{-1}$ are coarse mappings   then $f$ is called an
 {\it asymorphism}.
The coarse spaces
$(X, \mathcal{E})$, $(X^{\prime}, \mathcal{E}^{\prime})$
 are called {\it coarsely equivalent} if there exist large subsets
 $Y\subseteq X$, $Y^{\prime}\subseteq X$ such that
 $(Y, \mathcal{E}|_{Y})$
 and $(Y^{\prime}, \mathcal{E}^{\prime}|_{Y^{\prime}})$
  are asymorphic.

To conclude the coarse vocabulary, we take a family
$\{(X_{\alpha}, \mathcal{E}_{\alpha}) :  \alpha< \kappa\}$
  of coarse  spaces and define the
  {\it product}
 $P_{\alpha< \kappa}(X_{\alpha}, \mathcal{E}_{\alpha})$
  as the set $P_{\alpha< \kappa} X_{\alpha}$
 endowed with the coarse  structure with the base
 $P_{\alpha< \kappa} \mathcal{E}_{\alpha}$.
 If $\varepsilon_{\alpha}\in\mathcal{E}_{\alpha}$, $\alpha<\kappa$
  and  $x,y\in P_{\alpha<\kappa}X_{\alpha}$,
  $x=(x_{\alpha})_{\alpha<\kappa}$, $y=(y_{\alpha})_{\alpha<\kappa}$
  then $(x,y)\in (\varepsilon_{\alpha})_{\alpha<\kappa}$
  if and only if $(x_{\alpha}, y_{\alpha})\in\varepsilon_{\alpha}$
   for every $\alpha<\kappa$.

Let $\mathfrak{M}$ be a class of coarse spaces closed under asymorphisms.
We say that $\mathfrak{M}$  is a {\it variety} if $\mathfrak{M}$  is closed under formation of subspaces
$({\bf S} \mathfrak{M} \subseteq \mathfrak{M})$, coarse images $({\bf Q} \mathfrak{M} \subseteq \mathfrak{M})$
 and products $({\bf P} \mathfrak{M} \subseteq \mathfrak{M})$.

For an infinite cardinal $\kappa$, we say that a coarse space $(X, \mathcal{E})$
 is {\it $\kappa$-bounded} if every subset $Y\subseteq X$  such that $|Y|<\kappa$
  is bounded, and denote by $\mathfrak{M}_{\kappa}$
  the variety of all $\kappa$-bounded coarse spaces.

We denote by $\mathfrak{M}_{single}$ and $\mathfrak{M}_{bound} $ the variety of singletons and the variety of all bounded coarse spaces.

Then we get the  chain of varieties
$$\mathfrak{M}_{single}\subset \mathfrak{M}_{bound} \subset\ldots\subset \mathfrak{M}_{\kappa}\subset\ldots\subset\mathfrak{M}_{\omega}.$$

In  section 2, we prove that every variety of coarse  spaces lies in this chain and, in section 3, we  discuss some extensions   of this result to coarse spaces endowed with  additional algebraic structures.


\section{Results}

We recall that a family $\mathcal{I}$ of subsets of a set $X$ is an
{\it ideal} in the Boolean algebra $\mathcal{P}_{X}$ of all subsets of $X$ if
 $\mathcal{I}$ is closed under finite unions  and subsets. Every ideal $\mathcal{I}$
 defines the coarse structure with the base $\{\mathcal{E}_{A}: A\in\mathcal{I}\}$
  where $\mathcal{E}_{A}= (A\times  A)\cup \bigtriangleup_{X}$,
  so $B(x, \mathcal{E}_{A})=A$
  if $x\in A$ and $B(x, \mathcal{E}_{A})=\{x\}$
  if $x\in X\setminus A$.
We denote the obtained coarse space  by $(X, \mathcal{I})$.
For a cardinal $\kappa$, $[X]^{<\kappa}$ denotes the ideal
$\{Y\subseteq X: |Y|< \kappa\}$.
If $(X, \mathcal{E})$ is a coarse space,  the family $\mathcal{I}$  of all bounded
subsets of $X$ is an ideal. The coarse space $(X, \mathcal{I})$ is
called the {\it companion} of $(X, \mathcal{E})$.

Let $\mathcal{K}$  be a class of coarse spaces. We say that a
coarse space $(X, \mathcal{E})$ is {\it free} with respect to $\mathcal{K}$ if, for
 every  $(X^{\prime}, \mathcal{E}^{\prime})\in \mathcal{K}$ every mapping
   $f: (X, \mathcal{E})\longrightarrow (X^{\prime}, \mathcal{E}^{\prime})$
 is coarse. For example,
 $(X, [X]^{<\kappa})$
 is free with respect to the variety $\mathfrak{M}_{\kappa}$.
Since $(\kappa, [\kappa]^{<\kappa})\in  \mathfrak{M}_{\kappa}$
 but $(\kappa, [\kappa]^{<\kappa})\notin  \mathfrak{M}_{\kappa^{\prime}}$
 for each $\kappa^{\prime}> \kappa$,
  the inclusion $\mathfrak{M}_{\kappa^{\prime}}\subset \mathfrak{M}_{\kappa}$
  is strict.
  \vspace{5 mm}

{\bf  Lemma 1.}  {\it If a coarse space $(X,\mathcal{E})$  is free with  respect to a class
$\mathcal{K}$ then $(X,\mathcal{E})$
 is free with respect to ${\bf S} \mathcal{K}$, ${\bf Q} \mathcal{K}$, ${\bf P} \mathcal{K}$.
\vspace{4 mm}

Proof.}
We verify only the second statement.
Let  $(X^{\prime},\mathcal{E}^{\prime})\in\mathcal{K}$,
$(X^{\prime\prime},\mathcal{E}^{\prime\prime})\in {\bf Q}\mathcal{K}$,
 and
 $h: (X^{\prime},\mathcal{E}^{\prime})\longrightarrow (X^{\prime\prime},\mathcal{E}^{\prime\prime})$
  be a coarse surjective mapping.
   We take an arbitrary $f: X\longrightarrow X^{\prime\prime}$
   and choose $h^{\prime}: X\longrightarrow X^{\prime}$
   such that $f=hh^{\prime}$.
Since $(X,\mathcal{E})$
is free with respect to $\mathcal{K}$,
$h^{\prime}: (X,\mathcal{E})\longrightarrow (X^{\prime},\mathcal{E}^{\prime})$
 is coarse so $f$  is coarse as the  composition of the coarse mappings $h, h^{\prime}$.
 $ \ \ \  \Box$

\vspace{4 mm}

{\bf  Lemma  2.}  {\it Let $X$ be a set and let $\mathcal{K}$  be a class of coarse spaces, $\mathcal{K}\neq\mathfrak{M}_{single}$.
Then there exists a coarse structure $\mathcal{E}$ on $X$ such that $(X, \mathcal{E})\in {\bf SP} \mathcal{K}$ and $(X, \mathcal{E})$ is free with respect to $\mathcal{K}$.
\vspace{4 mm}

Proof.}
We take a set $S$ of all pairwise  non-asymorphic
coarse spaces $(X^{\prime}, \mathcal{E}^{\prime})\in \mathcal{K}$ such that
 $|X^{\prime}|\leq  |X|$  and enumerate all possible triplets
 $\{(X_{\alpha}, \mathcal{E}_{\alpha}, f_{\alpha}): \alpha <\lambda \}$ such that
 $(X_{\alpha}, \mathcal{E}_{\alpha})\in S$ and $f_{\alpha} : X\longrightarrow  X_{\alpha}$.
Then we consider the product
$P_{\alpha<\lambda}  (X_{\alpha}, \mathcal{E}_{\alpha})$ and define
$f: X \longrightarrow  P_{\alpha<\lambda} X_{\alpha}$   by
 $f(x)= (f _{\alpha}(x)) _{\alpha<\lambda}$.
Since $\mathcal{K} \neq \mathfrak{M }_{single}$,
 $f$ is injective so we can identify $X$ with $f(X) $   and consider the subspace
 $(X, \mathcal{E})$  of $P_{\alpha<\lambda}  (X_{\alpha}, \mathcal{E}_{\alpha})$ .
Clearly, $(X, \mathcal{E})\in {\bf SP} \mathcal{K}$.

To see that $(X, \mathcal{E})$  is free with respect to $\mathcal{K}$, it suffices
to verify that, for each  $(X^{\prime}, \mathcal{E}^{\prime})\in S$, every mapping
 $h: (X, \mathcal{E})\longrightarrow (X^{\prime}, \mathcal{E}^{\prime})$ is coarse. We take
 $\beta<\lambda$  such that  $(X^{\prime}, \mathcal{E}^{\prime})=(X_{\beta}, \mathcal{E}_{\beta}) $
  and $h=f_{\beta}$.
Then $f_{\beta}$ is the restriction to $X$ of the projection
$pr_{\beta}: P_{\alpha<\lambda} (X_{\alpha}, \mathcal{E}_{\alpha})\longrightarrow (X_{\beta}, \mathcal{E}_{\beta})$
  so $f_{\beta}$ is coarse.  $ \  \  \  \Box$
  \vspace{6 mm}

{\bf Theorem 1.}  {\it  For  every class $\mathcal{K}$ of coarse spaces, the smallest variety
Var  $\mathcal{K}$ containing $\mathcal{K}$ is ${\bf QSP}\mathcal{K}$.
\vspace{4 mm}

Proof.}The inclusion ${\bf QSP}\mathcal{K}\subseteq \mathcal{K}$ is evident. To prove
the inverse inclusion, we suppose that $\mathcal{K}\neq\mathfrak{M}_{single}$
 (this case is evident) and take an arbitrary   $(X^{\prime}, \mathcal{E}^{\prime})\in Var \ (\mathcal{K})$.
Then $(X^{\prime}, \mathcal{E}^{\prime})$ can be obtained from $\mathcal{K}$ by means of some finite sequence of operations ${\bf S, P, Q}$.
We use Lemma 2 to  choose a coarse space  $(X, \mathcal{E})\in {\bf SP}\mathcal{K}$, $|X|=|X^{\prime}|$ free with respect to $\mathcal{K}$. By Lemma 1, any bijection
$f: (X, \mathcal{E})\longrightarrow (X^{\prime}, \mathcal{E}^{\prime})$  is coarse so $(X^{\prime}, \mathcal{E}^{\prime})\in {\bf QSP}\mathcal{K}$. $ \  \  \  \Box$
\vspace{6 mm}

{\bf Theorem 2.}  {\it
Let $\mathfrak{M}$ be a variety of coarse spaces such that
$\mathfrak{M}\neq \mathfrak{M}_{single}$, $\mathfrak{M}\neq \mathfrak{M}_{bound}$.
Then there exists a cardinal $\kappa$  such that $\mathfrak{M}=\mathfrak{M}_{\kappa}$.

\vspace{4 mm}

Proof.} Since $\mathfrak{M}\neq \mathfrak{M}_{bound}$ and $\mathfrak{M}\neq \mathfrak{M}_{single}$, there exists the minimal cardinal $\kappa$ such that
$\mathfrak{M}$  contains an unbounded space of cardinality $\kappa$
 so $\mathfrak{M}\subseteq \mathfrak{M}_{\kappa}$.

To verify the inclusion $\mathfrak{M}_{\kappa}\subseteq \mathfrak{M}$,
we take a coarse space  $(X, \mathcal{E})\in\mathfrak{M}$ free with respect to $\mathfrak{M}$
and show that $(X, \mathcal{E})$
 is free with respect to  $\mathfrak{M}_{\kappa}$.
We prove that
$(X, \mathcal{E})=(X, [X]^{<\kappa})$.
If $|X|< \kappa$
then $(X, \mathcal{E})$
is bounded and the statement is evident.
Assume that $|X|\geq \kappa$
 but $(X, \mathcal{E})\neq(X, [X]^{<\kappa})$.
Assume that, for every $\varepsilon$, $\varepsilon=\varepsilon^{-1}$, the set
$S_{\varepsilon}= \{x\in X: |B(x,\varepsilon)|>1\}$ is bounded in $(X, \mathcal{E})$.
By the choice of $\kappa$, $|S_{\varepsilon}|< \kappa$ and $|B(x,\varepsilon)|=1$ for all
$x\in X\setminus S_{\varepsilon}$.  It follows that
$(X,\mathcal{E})=(X, [X]^{<\kappa})$.
Then there exists $\varepsilon\in\mathcal{E}$
such that the set
 $S_{\varepsilon}$
 is unbounded in  $(X, \mathcal{E})$.
We choose  a maximal by inclusion
 subset $Y\subset X$ such that
$B(y,\varepsilon)\cap B(y^{\prime},\varepsilon)=\emptyset$
for all distinct  $y, y^{\prime}\in Y$.
We observe that $Y$ is unbounded so
$|Y|\geq \kappa$.
We take an arbitrary $x_{0}\in X$  and choose a mapping $f: X \longrightarrow X$  such that
$f(y)=x_{0}$
 for each $y\in Y$ and $f $ is injective on $X\setminus Y$.  Since $(X, \mathcal{E}) $ is free with respect to
   $\mathfrak{M}$,
   the mapping
   $f: (X, \mathcal{E})\longrightarrow  (X, \mathcal{E}) $
    must be coarse. Hence, there exists
    $\varepsilon^{\prime}\in\mathcal{E}$
     such that
     $f(B(x,\varepsilon))\subseteq B(f(x),\varepsilon^{\prime})$
     for each $x\in X$.
It follows  that
$f(\cup_{y\in Y} B(y,\varepsilon))$
 is bounded in $(X, \mathcal{E})$.
We note that
$|f(\cup_{y\in Y} B(y,\varepsilon))|\geq\kappa$
 so $(X, \mathcal{E})$
  contains a bounded subset $Z$  such that  $|Z|=\kappa$.
Since $(X, \mathcal{E})$ is free with respect to $\mathfrak{M}$,  every
$(X^{\prime}, \varepsilon^{\prime})\in\mathfrak{M}$
 is a
 $\kappa^{+}$-bounded and we get a contradiction with the choice of $\kappa$.
 To conclude the proof, we take an arbitrary $(X, \mathcal{E}^{\prime})\in\mathfrak{M}_{\kappa}$
 and note that the identity mapping $id: (X, [X]^{<\kappa})\longrightarrow (X, \mathcal{E}^{\prime})$
 is coarse so $(X, \mathcal{E}^{\prime})\in\mathfrak{M}$.
 $ \  \  \  \Box$
\vspace{6 mm}

{\bf Remark 1.}
We note that $\mathfrak{M}_{single}$ is not closed under
coarse equivalence because each bounded  coarse space is coarsely equivalent to a singleton.
Clearly, $\mathfrak{M}_{bound}$  is closed
under coarsely equivalence. We show that the same is true for every variety
$\mathfrak{M}_{\kappa}$.
Let $(X, \mathcal{E})$
be a coarse space, $Y$ be a large subset of  $(X, \mathcal{E})$.
We assume that $(Y, \mathcal{E}|_{Y})\in\mathfrak{M}_{\kappa}$
 but  $(X, \mathcal{E})\notin\mathfrak{M}_{\kappa}$. Then $X$ contains an unbounded subset $Z$ such that $|Z|<\kappa$.
We choose $\varepsilon\in\mathcal{E}$
  such that $\varepsilon=\varepsilon^{-1}$
   and $X= B(Y, \mathcal{E})$.
For each $z\in Z$, we pick  $y_{z}\in Y$
such that  $z\in B(y_{z}, \mathcal{E})$.
We put  $Y^{\prime}= \{y_{z}\in Z\}$.
Since  $|Y^{\prime}|<\kappa$, $Y^{\prime}$ is bounded in $(Y, \mathcal{E}|_{Y})$.
It follows that $Z$  is bounded in
$(X, \mathcal{E})$, a contradiction with  the choice of $Z$.

We note also that every variety of coarse spaces is closed under formations of companions.
For $\mathfrak{M}_{single}$
and $\mathfrak{M}_{bound}$,
this is evident. Let $(X, \mathcal{E})\in\mathfrak{M}_{\kappa}$
and $\mathcal{I}$  is the ideal of all bounded subsets of $(X, \mathcal{E})$.
 Since $(X, [X]^{<\kappa})$
  is free with respect  to  $\mathfrak{M}_{\kappa}$, the identity mapping
  $id: (X, [X]^{<\kappa})\longrightarrow (X, \mathcal{E})$
    is coarse so $[X]^{<\kappa}\subseteq\mathcal{I}$
    and $(X, \mathcal{E})\in\mathfrak{M}_{\kappa}$.

\vspace{6 mm}

{\bf Remark 2.} Every metric $d$ on a set $X$  defines the coarse structure $\mathcal{E}_{d}$
 on $X$  with the base
 $\{(x,y): d(x,y)\leq n\} $,  $n\in\omega$.
A coarse structure $\mathcal{E}$  on $X$ is called {\it metrizable} if there exists a metric $d$ on  $X$
 such that  $\mathcal{E}=\mathcal{E}_{d}$.
By [9, Theorem 2.1.1], $\mathcal{E}$ is  metrizable if and only if $\mathcal{E}$  has a countable base.
From the coarse point of view, metric spaces are studding in    {\it Asymptotic Topology}, see [1].

We assume that a variety $\mathfrak{M}$ of coarse space contains an unbounded metric space $(X, d)$ and show that
$\mathfrak{M}=\mathfrak{M} _{\omega}$.
We choose a countable unbounded subset $Y$ of $X$ and note that $(Y,d)\notin \mathfrak{M}_{\kappa}$
for $\kappa>\omega$
 so $(Y,d)\in \mathfrak{M} _{\omega}\setminus \mathfrak{M} _{\kappa}$
  and the variety generated by $(X,d)$  is $\mathfrak{M} _{\omega}$.


\section{Comments}

1.	Let $G$  be a group with the identity $e$.
An ideal $\mathcal{I}$  in $\mathcal{P}_{G}$  is called a {\it group ideal} if $[G]^{< \omega}\subseteq \mathcal{I}$  and $AB ^{-1} \in \mathcal{I}$  for all $A, B \in \mathcal{I}$.

Let $X$  be a $G$-space with the action $G\times X\longrightarrow  X$,  $(g,x)\longmapsto gx$.
We assume that $G$ acts on $X$ transitively, take a group ideal $\mathcal{I}$ on $G$
and consider the coarse structure
 $\mathcal{E}(G, \mathcal{I}, X)$ on $X$ with the base  $\{\varepsilon _{A}: A\in \mathcal{I},  e\in A\}$,
 $\varepsilon _{A}= \{(x, gx): x\in X, g\in  A\}$.
Then $B(x, \varepsilon _{A})= Ax$, $Ax= \{gx: g\in A\}$.

By [5, Theorem 1], for every coarse structure $\mathcal{E}$  on $X$,
there exist a group $G$ of permutations of $X$  and a group ideal $\mathcal{I}$
 in $\mathcal{P}_{G}$  such that $\mathcal{E}=\mathcal{E}(G, \mathcal{I}, X)$.

Now let $X=G$ and $G$ acts on $X$ by the left shifts
$x\longmapsto gx$, $g\in G$. We denote
$(G, \mathcal{E}(G, \mathcal{I}, G))$  by $(G, \mathcal{I})$
and say that $(G,\mathcal{I})$  is a {\it right coarse group}.
 If $\mathcal{I}=[G]^{<\omega}$  then $(G, \mathcal{I})$ is called  {\it a finitary right coarse group}.
In the metric form, these structures on finitely  generated groups play an important role in  {\it Geometric Group Theory}, see [2, Chapter 4].

A group $G$ endowed with a coarse structure $\mathcal{E}$  is a right coarse group if and only if, for every $\varepsilon\in \mathcal{E}$,  there exists $\varepsilon^{\prime}\in \mathcal{E}$ such that
 $(B(x,\varepsilon))g \subseteq B(xg, \varepsilon^{\prime})$ for all $x,g\in G$.
For group ideals and coarse structures on groups see \cite{b8} or [9, Chapter 6].

\vspace{6 mm}

2.	A class  $\mathfrak{M}$  of right coarse groups is called a {\it variety}
if $\mathfrak{M}$ is closed under formation of subgroups, coarse homomorphic images and products.

Let $\mathcal{K}$ be a class of right coarse groups,  $G$  be a group generated by a subset $X\subset G$. We say that a right coarse group $(G, \mathcal{I})$ is {\it free} with respect to  $\mathcal{K}$  if, for every
$(G^{\prime}, \mathcal{I}^{\prime})\in\mathcal{K}$,
any mapping  $X\longrightarrow G^{\prime}$  extends to the coarse homomorphism
$(G, \mathcal{I})\longrightarrow (G^{\prime}, \mathcal{I}^{\prime})$.
Then Lemma 1, Lemma 2 and Theorem 1 hold for right coarse groups in place of coarse spaces.

Let $\mathfrak{M}$  be a variety of right coarse groups. We take
an arbitrary $(G, \mathcal{I})\in \mathfrak{M}$, delete the coarse structure
 on $G$  and the class $\mathfrak{M}^{\flat}$ of groups. If $(G, \mathcal{I})\in \mathfrak{M}$ then
 $(G, \mathcal{P}_{G})\in \mathfrak{M}$.
It follows that $\mathfrak{M}^{\flat}$ is  a variety of groups.

Now let  $\mathcal{G} $  be a variety of group different from the variety of singletons. We denote by
$\mathcal{G} _{bound}$ the   variety of right coarse groups
$(G, \mathcal{P}_{G}),  G\in \mathcal{G}$.
For an infinite cardinal   $\kappa$, we denote by
$\mathcal{G} _{\kappa}$  the variety of all $\kappa$-bounded right coarse groups  $(G,\mathcal{I})$, $G\in \mathcal{G}$.

Let $\mathfrak{M}$ be a variety of right coarse groups such that $\mathfrak{M}^{\flat}= \mathcal{G}$.
In contrast to Theorem 2, we can not state  that
 $\mathfrak{M}$  lies in the chain
$$\mathcal{G} _{bound} \subset\ldots\subset \mathcal{G} _{\kappa}\subset\ldots\subset \mathcal{G} _{\omega}.$$
If $G$  is a group  of cardinality $\kappa$  and $G\in \mathcal{G} $ then
$(G, [G]^{<\kappa})\in \mathcal{G}_{\kappa}\setminus \mathcal{G}_{\kappa^{\prime}} $
for each  $\kappa^{\prime}>\kappa$.
Hence, all inclusions in above chain are strict.

Let $\mathcal{G}$ be variety of all groups. We define the variety $\mathfrak{M}$
 of right coarse groups by the rule: $(G, \mathcal{I})\in\mathfrak{M}$ if and only if  $\mathcal{I}$  contains all finitely generated subgroups of $G$.
Then   $\mathfrak{M}^{\flat}= \mathcal{G}$  but $\mathfrak{M}$  does not lie in above chain.

\vspace{6 mm}

3. let  $\Omega$ be a signature, $A$ be an  $\Omega$-algebra and
$\mathcal{E}$ be a coarse structure on $A$.
We say that $A$ is a {\it coarse $\Omega$-algebra}  if every    $n$-ary operation from $\Omega$ is coarse as the mapping
$(A,\mathcal{E})^{n}\longrightarrow (A,\mathcal{E})$.
 We note that each coarse group is a right coarse group but the converse statement needs not to be true, see [9, Section 6.1].

A class $\mathfrak{M}$  of coarse $\Omega$-algebra is called a {\it variety}  if $\mathfrak{M}$  is closed under formation of subalgebras, coarse homomorphic images and  products. Given a variety  $\mathfrak{M}$
of coarse algebras, the class $\mathfrak{M}^{\flat}$
 of all $\Omega$-algebras $A$ such that $(A, \mathcal{E})\in \mathfrak{M}$  is a variety of  $\Omega$-algebras.

Let $A$ be a variety of  $ \  \Omega$-algebras  different from the
variety of singletons. We denote by $\mathcal{A} _{bound}$ the variety of coarse algebras $(A, \mathcal{P }_{A})$,
$A \in \mathcal{A}$.
For an infinite cardinal, $\kappa$  we denote by $\mathcal{A}_{\kappa}$  the variety of all
$\kappa$-bounded $\Omega$-algebras $(A, \mathcal{E})$   such that   $A \in \mathcal{A}$, and get the chain.
$$\mathcal{A} _{bound} \subseteq \ldots  \subseteq  \mathcal{A}_{\kappa} \subseteq  \ldots  \subseteq  \mathcal{A}_{\omega} ,$$  but,
 we can not state that all inclusions are strict.
In the case of course groups, this is so because
 each non-trivial variety of groups contains some Abelian group $A$ of cardinality $\kappa$
 and the coarse group $(A, [A] ^{<\kappa)}$ is $\kappa$-bounded but not $\kappa^{+}$-bounded.

\vspace{6 mm}

4. A class $\mathfrak{M}$ of topological $\Omega$-algebras (with regular topologies)
 is called a {\it variety} ({\it a wide variety}) if $\mathfrak{M}$ is closed under formation of  closed subalgebras   (arbitrary subalgebras), continuous homomorphic images and products.
  The wide varieties and varieties are characterized  syntactically by the limit laws \cite{b11}   and  filters \cite{b6}.
In our coarse case, the part of filters play the ideals $[X] ^{<\kappa}$.

There are only two wide varieties of topological spaces, the variety of   singletons and the variety of all topological spaces, but there is  a plenty of  varieties of topological spaces. The variety of coarse  spaces $\mathfrak{M}_{\kappa}$ is a twin of the varieties of topological spaces in which every subset of cardinality $<\kappa$  is compact.
We note also that $\mathcal{G}_{\kappa} $
 might be considered as a  counterpart of the variety $T( \kappa ) $ of topological groups from \cite{b3},
 $ G\in T( \kappa)$ if and only if each neighborhood of $e$  contains  a normal subgroup of index  strictly  less then $\kappa$.
\vspace{6 mm}

5. A class  $\mathfrak{M}$  of uniform spaces is called a {\it variety} if $\mathfrak{M}$ is closed under formation of  subspaces, products and uniformly continuous images.
For an  infinite cardinal $\kappa$, a uniform space $X$ is called  $\kappa$-bounded if $X$ can be covered by $< \kappa$ balls of arbitrary small radius. Every variety of uniform spaces different from varieties of singletons and all spaces  coincides with the variety of $\kappa$-bounded spaces for some $\kappa$, see \cite{b4}.
I thank Miroslav Hu$\check{s}$ek for this reference.

\vspace{5 mm}

6.  On varieties of bornological spaces. An ideal $\mathcal{I}$ in $\mathcal{P}_{X}$
 is called a {\it bornology  } if  $\bigcup \mathcal{I}=X$.
A set $X$, endowed with a bornology $\mathcal{I}$  is called a {\it bornological space}.
Each $A\in\mathcal{I}$ is called {\it bounded},  $(X,\mathcal{I})$
is bounded if $X\in\mathcal{I}$.
For an infinite cardinal $\kappa$, $(X,\mathcal{I})$  is called  $\kappa$-{\it bounded} if
$[X]^{<\kappa}\subseteq\mathcal{I}$.

A mapping $f: (X,\mathcal{I})\longrightarrow (X^{\prime},\mathcal{I}^{\prime})$
 is called {\it bornologous  } if $f(A)\in \mathcal{I}^{\prime}$
  for each  $A\in\mathcal{I}$.

A class of bornological space closed under subspaces,
 products and bornologous images is called a {\it variety}.
Repeating the first part of the proof of Theorem 2,
 we conclude that each variety of bornological  spaces is either variety of singletons,
 or variety of all bounded spaces or variety of all $\kappa$-bounded spaces for some infinite cardinal $\kappa$.


\vskip 5pt

CONTACT INFORMATION

I.~Protasov: \\
Faculty of Computer Science and Cybernetics  \\
        Kyiv University  \\
         Academic Glushkov pr. 4d  \\
         03680 Kyiv, Ukraine \\ i.v.protasov@gmail.com

\end{document}